%% file: machine_learning_weak_approx_arXiv_revised.tex
\documentclass[a4paper,1	2pt]{article}

\makeatletter
 
  \@addtoreset{equation}{section}
 \makeatother
\usepackage{amsfonts}
\usepackage{amssymb}
\usepackage{mathrsfs}
\usepackage{amsmath,amsthm}
\usepackage{amsmath}
\usepackage{comment}
\usepackage{algorithm}
\usepackage{algorithmic}
\usepackage{shuffle}

\usepackage[dvipdfmx]{graphicx}
\usepackage{color}
\usepackage{indentfirst}
\usepackage{here}

\begin{document}
\input{tymacro_3}

\newtheorem{thm}{Theorem}
\newtheorem{lemma}{Lemma}
\newtheorem{prop}{Proposition}
\newtheorem{defn}{Definition}
\newtheorem{rem}{Remark}
\newtheorem{step}{Step}
\newtheorem{cor}{Corollary}
\newtheorem{assump}{Assumption}

\newcommand{\Cov}{\mathop {\rm Cov}}
\newcommand{\Var}{\mathop {\rm Var}}
\renewcommand{\E}{\mathop {\rm E}}
\newcommand{\const }{\mathop {\rm const }}
\everymath {\displaystyle}

\newcommand{\ruby}[2]{
\leavevmode
\setbox0=\hbox{#1}
\setbox1=\hbox{\tiny #2}
\ifdim\wd0>\wd1 \dimen0=\wd0 \else \dimen0=\wd1 \fi
\hbox{
\kanjiskip=0pt plus 2fil
\xkanjiskip=0pt plus 2fil
\vbox{
\hbox to \dimen0{
\small \hfil#2\hfil}
\nointerlineskip
\hbox to \dimen0{\mathstrut\hfil#1\hfil}}}}

\def\qedsymbol{$\blacksquare$}
\renewcommand{\thefootnote }{\fnsymbol{footnote}}
\renewcommand{\refname }{References}

\everymath {\displaystyle}

\title{
A machine learning solver for high-dimensional integrals: 
 Solving Kolmogorov PDEs by stochastic weighted minimization and stochastic gradient descent through a high-order weak approximation scheme of SDEs with Malliavin weights
} 
\author{Riu Naito\footnote{Asset Management One, Tokyo, Japan} \ and \ Toshihiro Yamada\footnote{Hitotsubashi University, Tokyo, Japan} \footnote{Japan Science and Technology Agency (JST), Tokyo, Japan}}
\date{First version: December 22, 2020, \ This version: February 11, 2021} 
\maketitle

\abstract 
The paper introduces a very simple and fast computation method for high-dimensional integrals 
to solve high-dimensional Kolmogorov partial differential equations (PDEs). The new machine learning-based method is obtained by solving a stochastic weighted minimization with stochastic gradient descent which is inspired by a high-order weak approximation scheme for stochastic differential equations (SDEs) with Malliavin weights. Then solutions to high-dimensional Kolmogorov PDEs or expectations of functionals of solutions to high-dimensional SDEs are accurately approximated without suffering from the curse of dimensionality. Numerical examples for PDEs and SDEs up to 100 dimensions are shown by using  second and third-order discretization schemes in order to demonstrate the effectiveness of our method. 

\section{Introduction}
In the paper, we introduce a novel machine learning-based computation method for high-dimensional integrals to solve high-dimensional Kolmogorov partial differential equations (PDEs) (\cite{F}\cite{KS}\cite{Taylor}\cite{BBGJJ}) without suffering from ``the curse of dimensionality" (\cite{TrWe}). 
 
 Let $T>0$ and $d \in \mathbb{N}$  be an integer which is assumed to be high. Let $u:[0,T] \times \mathbb{R}^d \to \mathbb{R}$ be the solution to the following Kolmogorov partial differential equation (PDE):  
\begin{eqnarray}
&\partial_t u(t,x)+\langle b(x), \nabla_x u(t,x) \rangle_{\mathbb{R}^d}+\frac{1}{2}\mathrm{Trace}_{\mathbb{R}^d}[\sigma\sigma^\top(x) \nabla^2_x u(t,x)]=0, \ \ t<T, \  x \in \mathbb{R}^d, \ \ \ \\
& u(T,x)=f(x), \ \  x \in \mathbb{R}^d, \label{PDE}
\end{eqnarray}
where $b$ and $\sigma$ are some appropriate functions on $\mathbb{R}^d$ valued in $\mathbb{R}^d$ and $\mathbb{R}^{d \otimes d}$, respectively, and $f$ is a function on $\mathbb{R}^d$ valued in $\mathbb{R}$. The purpose of this paper is to show a new computation method of 
\begin{eqnarray}
u(0,x_0)=\int \cdots \int_{\mathbb{R}^d} f(y) p(T,x_0,y)dy_1\cdots dy_d,  \label{target_value_pde}
\end{eqnarray} 
for $x_0 \in \mathbb{R}^d$, where the function $(t,x,y)\mapsto p(T-t,x,y)$ is the fundamental solution (with the condition $\textstyle{\lim_{t \to T}p(T-t,x,\cdot)=\delta_{x}(\cdot)}$ where $\delta_{x}$ is the Dirac mass at the point $x$). The computation problem of $u(0,x_0)$ is reformulated using language of probability theory. Let $\{ W_t \}_{t\geq 0}$ be a $d$-dimensional Brownian motion on a probability space and consider the solution of the  following SDE: 
\begin{eqnarray}
dX_t^x=b(X_t^x)dt+\sum_{i=1}^d \sigma_i(X_t^x) dW_t^i, \ \ X_0^x=x \in \mathbb{R}^d.
\end{eqnarray}
Then, it holds that 
\begin{eqnarray}
u(0,x_0)=E[f(X_T^{x_0})]=\int \cdots \int_{\mathbb{R}^d} f(y) p(T,x_0,y)dy_1\cdots dy_d
\end{eqnarray}
through the Feynman-Kac formula. 

As PDE solvers, finite element method and finite difference method (\cite{Thomee}\cite{LT}\cite{Braess}) are known as classical schemes. However, these methods basically work in $d \leq 2$, $3$ but fail if the dimension $d$ is high since the cost grows exponentially. The use of the Monte Carlo method (\cite{MU}\cite{HH}\cite{Cafl}\cite{G}) with the Euler-Maruyama scheme (\cite{Maru}) (the first-order discretization of stochastic differential equations (SDEs)) is a possible approach to solve $u(0,x_0)$. Although the Monte Carlo approximation is regarded as the standard method today, its computation time is not fast as usual. Thus, an efficient numerical integration method and an efficient time-discretization method are required to solve $u(0,x_0)$. As an accurate numerical method for multidimensional integrals, the quasi Monte Carlo (QMC) method is known (\cite{H}\cite{So}\cite{Faure}\cite{Nied}\cite{Owen}\cite{DKS}). However, the accuracy of QMC becomes worse when the dimension $d$ becomes high in general. In other words, QMC works in relatively low-dimensional cases. In terms of time-discretization, employing a higher order method will be a refined approach. On the other hand, higher order time-discretization methods have essential difficulty in the implementation due to ``L\'evy areas" or higher order iterated stochastic integrals in stochastic Taylor expansion in high-dimensional cases (\cite{Wik}\cite{KP}). Also, for higher order methods, the regularity of the terminal condition $f$ should be carefully discussed so that stochastic Taylor expansion is appropriately applied. Then, constructing a numerical integration method and a time-discretization method is still important or is of interest to solve $u(0,x_0)$ in high-dimensional cases. 

Based on the background mentioned above, we introduce a new machine learning-based computation method inspired by the pioneering works in \cite{BBGJJ}\cite{BBCJN}, combining with high-order weak approximation with Malliavin weights proposed in \cite{Y}\cite{NY}\cite{IY}. The high-order weak approximation method (extension of \cite{TY}) based on Malliavin calculus is applied for solving high-dimensional SDEs. 
The Malliavin weight plays a role in the method, which is constructed with Malliavin's integration by parts (\cite{IW}\cite{M}\cite{N}). The main advantages of the use of the high-order method are the followings: 
\begin{enumerate}
\item the method provides ``any order" discretization algorithm;
\item different from other higher order time-discretization methods, the required number of random variables 
 in each time step in the implementation is always $d$ (=the same as the dimension of the underlying SDE as in the Euler-Maruyama scheme); 
\item the method can be applied to smooth and nonsmooth functionals of solutions to SDEs 
\end{enumerate}

These will be distinguishing features of the method comparing with other discretization methods. Indeed, even if the diffusion matrix does not satisfy the commutativity condition (i.e. even for the case that the L\'evy areas appear), the method provides universal algorithm for general order $m$, although other higher order methods need a specific implementation scheme for a specific weak order $m$ (or have no implementation scheme). Also, even if their implementation schemes is obtained, they require many random variables much more than $n \times d$ where $n$ is the number of discretization, which suggests that their computational costs become high. However, the Malliavin-based weak approximation need only $n \times d$-random variables in the implementation. This is because the Malliavin weight is explicitly given by certain polynomials of $d$-dimensional Brownian motion. Furthermore, the Malliavin-based weak approximation adopts very small $n$ since it provides higher order discretization scheme. 

The paper propose a new computation scheme based on the Malliavin-based weak approximation and stochastic gradient descent in order to obtain an extension of the methods in \cite{BBGJJ}\cite{BBCJN}. 
We found that stochastic gradient descent method is suitable for implementation of the Malliavin-based weak approximation in \cite{IY}\cite{NY}\cite{Y}. The high-dimensional integrals appearing  in weak approximation of high-dimensional SDEs or high-dimensional Kolmogorov PDEs are easily solved with good accuracy through ``Malliavin weighted" minimization. We check that the proposed method works well in high-dimension case such as $d=100$. 

After reviewing the weak approximation with Malliavin weights (section 2), we provide a new stochastic gradient descent-based method (section 3). Numerical examples are shown in section 4. Section 5 concludes the method and mentions future works.  

\section{Weak approximation of SDE}
We formulate our setting. Let $(\Omega,{\cal F},P)$ be a $d$-dimensional Wiener space, i.e. $\Omega=C_{(0)}([0,T]; \mathbb{R}^d)$$=\{ \omega:[0,T]\to \mathbb{R}^d; \ \omega(0)=0, \ \omega \ \mbox{is \ continuous} \}$, ${\cal F}$ is the Borel field over $\Omega$ and $P$ is the Wiener measure.  
Let $W=\{ W_t \}_{t\geq 0}$ be a $d$-dimensional Brownian motion on the Wiener space and consider the solution of the following SDE: 
\begin{eqnarray}
X_t^x=x+\int_0^t b(X_s^x)ds+\sum_{i=1}^d \int_0^t \sigma_i (X_s^x) dW_s^i,
\end{eqnarray}
where $b,\sigma_i \in C_b^\infty(\mathbb{R}^d;\mathbb{R}^d)$, $i=1,\cdots,d$. We assume that there is $\varepsilon_0>0$ such that $\sigma(x)\sigma(x)^{\top} \geq \varepsilon_0 I_{\mathbb{R}^{d \otimes d}}$ for all $x \in \mathbb{R}^d$.

Let $\{ P_t \}_{t\geq 0}$ be a semigroup of linear operators given by 
\begin{eqnarray}
(P_t f)(x)=E[f(X_t^x)], \ \ \ t\geq 0, \ x \in \mathbb{R}^d, 
\end{eqnarray}
for a Lipschitz continuous function $f:\mathbb{R}^d \to \mathbb{R}$. The purpose of the paper is to give a new computation scheme for the value:  
\begin{eqnarray}
(P_T f)(x_0)=E[f(X_T^{x_0})]
\end{eqnarray}
with $T>0$, $x_0 \in \mathbb{R}^d$ and a Lipschitz continuous function $f:\mathbb{R}^d \to \mathbb{R}$, which corresponds to (\ref{target_value_pde}). See \cite{F}\cite{KS}\cite{MT}\cite{Taylor}\cite{ELV} for stochastic calculus, PDE and related topics for instance.  

First, we recall the following approximation scheme (\cite{Y}\cite{NY}\cite{IY}).
\begin{thm}
Let $m\geq 1$ and $\{ Q_{t}^{(m)} \}_{t>0}$ be the linear operators which has the form $(Q_{t}^{(m)} f)(x)=E[f(\bar{X}_t^x) M^{(m)}_t(x,W_t) ]$ for $t>0$, $x \in \mathbb{R}^d$, and a Lipschitz continuous function $f: \mathbb{R}^d \to \mathbb{R}$, with the one-step Euler-Maruyama approximation $$\textstyle{\bar{X}_t^x=x+b(x)t+\sum_{i=1}^d \sigma_i(x)W^i_t}$$ and the Malliavin weight $M^{(m)}_t(x,W_t)$ such that  
\begin{eqnarray}
\Big\| P_t f - Q_{t}^{(m)} f \Big\|_{\infty}=O(t^{m+1}), 
\end{eqnarray} 
for all $f \in C_b^\infty(\mathbb{R}^d)$. Then, for $T\geq 1$ and $m\geq 1$, we have 
\begin{eqnarray}
\Big\| P_T f - (Q^{(m)}_{T/n})^n f \Big\|_{\infty} = O(n^{-m}), \label{WA_thm}
\end{eqnarray}
for any Lipschitz continuous function $f: \mathbb{R}^d \to \mathbb{R}$. 
\end{thm}

\begin{rem}
The approximation (\ref{WA_thm}) also holds for any bounded measurable function $f: \mathbb{R}^d \to \mathbb{R}$.
\end{rem}
 The construction of the weak approximation is based on Malliavin calculus, in particular, the error estimate of the discretization is obtained by employing the theory of Kusuoka-Stroock \cite{KS}. The Malliavin weights  are explicitly calculated using Malliavin's integration by parts on the Wiener space (\cite{IW}\cite{M}\cite{N}). For example, 
for $m=2$, 
 the Malliavin weight $M_t^{(2)}(x,W_t) $ is simply given as follows: 
\begin{align}
& M_t^{(2)}(x,W_t) 
=1+\sum_{i_1,i_2=0}^{d} \sum_{i_3,i_4=1}^{d} \frac{1}{2t} L_{i_1}V^{i_4}_{i_2}(x)
(\sigma^{-1})^{i_3i_4}(x)\nn\\
& \ \ 
 \{ W_t^{i_1}W_t^{i_2}W_t^{i_3}-tW_t^{i_3}{\bf 1}_{\{i_1=i_2\neq 0 \}}-tW_t^{i_1}{\bf 1}_{\{i_2=i_3\neq 0 \}}-tW_t^{i_2}{\bf 1}_{\{i_1=i_3 \neq 0 \}} \}\nn\\
& \ \  
 + \sum_{i_1,\cdots,i_6=1}^{d}  \frac{1}{4} L_{i_1}V^{i_4}_{i_2}
(x)L_{i_1}V^{i_6}_{i_2}(x) (\sigma^{-1})^{i_3 i_4}(x)(\sigma^{-1})^{i_5i_6}(x)
  \{W_t^{i_3}W_t^{i_5}-t{\bf 1}_{\{i_3=i_5  \}} \},  \label{M_weight}
\end{align}
where $V_0(\cdot)=b(\cdot)$, $V_i(\cdot)=\sigma_i(\cdot)$, $i=1,\cdots,d$, and 
$L_i$, $i = 0, 1,\cdots, d$ are the differential operators appearing in the It\^o Taylor expansion (\cite{KP}): for $\varphi \in C_b^{\infty}(\mathbb{R}^d)$, $x \in \mathbb{R}^d$, 
\begin{eqnarray}
L_0 \varphi(x)&=&\sum_{k=1}^d b^k(x) \frac{\partial \varphi}{\partial x_k}(x)+\frac{1}{2} \sum_{k,\ell,j=1}^d \sigma_j^k(x)\sigma_j^\ell(x) \frac{\partial^2 \varphi}{\partial x_k \partial x_\ell}(x),\\
L_i \varphi(x)&=&\sum_{k=1}^d \sigma_i^k(x) \frac{\partial \varphi}{\partial x_k}(x), \ \ \ i=1,\cdots,d. 
\end{eqnarray}
More higher order cases are also obtained by adding some polynomials of Brownian motions $W^i$, $i=1,\cdots,d$. Note that the case $m=1$ corresponds to the Euler-Maruyama scheme (\cite{Maru}) (for nonsmooth test functions), i.e. $M_t^{(1)}(x,W_t)\equiv 1$, which can be checked in (\ref{basic_formula}) later.

In the next section, we show how the numerical value $(Q^{(m)}_{T/n})^n f(x_0)$ of the operator splitting approximation can be simply solved by using machine learning. 

\section{Machine learning-based method}
Hereafter, let $f:\mathbb{R}^d \to \mathbb{R}$ be a Lipschitz continuous function. We give a machine learning-based method for computing $(P_T f)(x_0)=E[f(X_T^{x_0})]$ for a fixed $x_0 \in \mathbb{R}^d$. 
We first note that the approximation $(Q^{(m)}_{T/n})^n f(x_0)$ can be explicitly represented using the Euler-Maruyama discretization and the Malliavin weight. 

Let $\{\bar{X}^{(n)}_{kT/n}\}_{0\leq k \leq n}$ be the Euler-Maruyama discretization: for $1\leq k \leq n$, 
\begin{align}
\bar{X}^{(n)}_{kT/n}=\bar{X}^{(n)}_{(k-1)T/n}+b(\bar{X}^{(n)}_{(k-1)T/n})T/n+\sum_{i=1}^d \sigma_i(\bar{X}^{(n)}_{(k-1)T/n}) \{ W^i_{kT/n}-W^i_{(k-1)T/n} \},
\end{align}
with $\bar{X}^{(n)}_{0}=x_0$. Then, it holds that 
\begin{eqnarray}
(Q^{(m)}_{T/n})^n f(x_0) = E[ f (\bar{X}^{(n)}_{T}) \prod_{i=1}^{n} M^{(m)}_{T/n}(\bar{X}^{(n)}_{(i-1)T/n},W_{iT/n}-W_{(i-1)T/n})  ]. \label{basic_formula}
\end{eqnarray}
We have the following new representation with the ``Malliavin weighted" minimization. 

\begin{prop}
It holds that
\begin{align}
 &(Q^{(m)}_{T/n})^n f(x_0) \nn\\
=&  \mathrm{argmin}_{v\in \mathbb{R}} E[ | v - f (\bar{{X}}^{(n)}_{T}) \prod_{i=1}^{n} M^{(m)}_{T/n}(\bar{{X}}^{(n)}_{(i-1)T/n},W_{iT/n}-W_{(i-1)T/n})  |^2 ].
\end{align}
\end{prop}
\noindent
{\it Proof}. 
First, note that $\Big\| \prod_{i=1}^{n} M^{(m)}_{T/n}(\bar{{X}}^{(n)}_{(i-1)T/n},W_{iT/n}-W_{(i-1)T/n}) \Big\|_2 \leq C_{b,\sigma,T}$ for all $n\geq 1$ (see \cite{IY} for instance). Then, by the property of expectation (see Lemma 2.1 of Beck et al. (2018) \cite{BBGJJ} for example), we have 
\begin{eqnarray}
&&
E[ | (Q^{(m)}_{T/n})^n f(x_0)  - f (\bar{{X}}^{(n)}_{T}) \prod_{i=1}^{n} M^{(m)}_{T/n}(\bar{{X}}^{(n)}_{(i-1)T/n},W_{iT/n}-W_{(i-1)T/n})  |^2 ] \nn\\
&=&
\inf_{v \in \mathbb{R}} E[ | v  - f (\bar{{X}}^{(n)}_{T}) \prod_{i=1}^{n} M^{(m)}_{T/n}(\bar{{X}}^{(n)}_{(i-1)T/n},W_{iT/n}-W_{(i-1)T/n})  |^2 ]. \ \ \ \Box
\end{eqnarray}

Using the above representation, the target $P_T f(x_0)$ is approximated by 
the (plain-vanilla) stochastic gradient descent as follows. Let $J \in \mathbb{N}$ be 
the number of iteration steps. By simulating $J$ independent trajectories of $W$ (and then of $\bar{{X}}^{(n)}$), for $j=1,\cdots,J$, we define
\begin{align}
\psi_j(\theta)= \Big| \theta  - f (\bar{{X}}^{(n)}_{T}(\omega_j)) \prod_{i=1}^{n} M^{(m)}_{T/n}(\bar{{X}}^{(n)}_{(i-1)T/n}(\omega_j),W_{iT/n}(\omega_j)-W_{(i-1)T/n}(\omega_j))  \Big|^2
\end{align}
and compute 
\begin{eqnarray}
\vartheta_j^\ast &=& \vartheta_{j-1}^\ast - \gamma_{j} \cdot (\nabla_\theta \psi_j)( 
\vartheta_{j-1}^\ast), \ \  \gamma_j \in (0,1). 
\end{eqnarray}
Then, the following approximation holds: 
\begin{eqnarray}
P_T f(x_0) &\approx& \vartheta_J^\ast. 
\end{eqnarray}
The above procedure is refined with mini-batches as well as {\it Adam} with mini-batches (see \cite{KB}\cite{R}) which is used in the numerical experiments.\\
 
\section{Algorithm \& Numerical examples}
We list the algorithm of the machine learning-based scheme used in the section where the Adam with mini-batches is employed to enhance the method in the previous section, which can be simply implemented with Python by partially making use of the code in \cite{BBGJJ} or \cite{BBCJN}. \\

\begin{algorithm}[H]                      
\caption{Machine learning-based high-order weak approximation method}\label{alg}                          
{{\begin{algorithmic}                  
\REQUIRE $M \in \mathbb{N}$ (batch size), $J\in \mathbb{N}$ (number of train-steps) $n \in \mathbb{N}$ (number of time-steps), $\gamma \in (0,1)$ (learning rate) 
\STATE Define the Malliavin weight function $(\xi,\Delta W) \mapsto M_{T/n}^{(m)}(\xi,\Delta W)$ of $m$-th order weak approximation
    \FOR {$j=1$ to $J$}
	\FOR {$\ell=1$ to $M$}
		\STATE $\bar{X}_{0}^{(n),\ell,j}=x$
		\FOR {$i=0$ to $n-1$}
			\STATE  $\Delta W^{\ell,j}_{(i+1)T/n}=(\Delta W^{1,\ell,j}_{(i+1)T/n},\cdots,
			\Delta 
			W^{d,\ell,j}_{(i+1)T/n})$, i.i.d $\Delta W^{k,\ell,j}_{(i+1)T/n} \sim N(0,T/n)$, 
			$k=1,\cdots,d$, 
			
			\STATE $\bar{X}_{(i+1)T/n}^{(n),\ell,j}=\bar{X}_{iT/n}^{(n),\ell,j}+
							b(\bar{X}_{iT/n}^{(n),\ell,j})T/n+\sum_{h=1}^d \sigma_k(\bar{X}_{iT/n}
							^{(n),\ell,j})\Delta W^{k,\ell,j}_{(i+1)T/n}$
			\STATE Compute Malliavin weight $M_{T/n}^{(m)}(\bar{X}_{iT/n}^{(n),\ell,j},\Delta W^{\ell,j}_{(i+1)T/n})$
		\ENDFOR
	\ENDFOR 
	\STATE Loss $\psi_j(\theta)=\frac{1}{M}\sum_{\ell=1}^M 
	\{ \theta- f(\bar{X}_{T}^{(n),\ell,j}) \prod_{i=1}^n M_{T/n}^{(m)}(\bar{X}_{(i-1)T/n}^{(n),\ell,j},\Delta W_{iT/n}^{\ell,j})\}
	^2$
	\STATE Update the optimized parameter $\vartheta_j^\ast$ using stochastic gradient descent 
    \ENDFOR
\STATE Return $\vartheta_J^\ast$
\end{algorithmic}}}
\end{algorithm}
\newpage

 We illustrate effectiveness of the machine learning-based high-order weak approximation through numerical examples for high-dimensional PDEs. We consider the following Kolmogorov PDE inspired by a high-dimensional model appearing in financial mathematics: 
\begin{eqnarray*}
\partial_t u(t,x)+\frac{1}{2}\sum_{i=1}^{d}\sigma_i^2 x_i^2 \partial_{x_i}^2 u(t,x)=0, \ \ \ \ u(T,x)=f(x),
\end{eqnarray*}
for a function $f:\mathbb{R}^d \to \mathbb{R}$ and estimate $u(0,x)$ of the solution $u:[0,T] \times \mathbb{R}^d \to \mathbb{R}$. The dimension $d$ is assumed to be up to $100$. All numerical experiments of this 
section are performed in Python using Tensorflow 1.12 from Google Colab.\\
  
\subsection{10-dimension case}  
We first check the numerical results in $10$-dimension case. Let $d=10$, $\sigma=0.2$, 
$T=2.0$, $f=f_K:\mathbb{R}^d \to \mathbb{R}$ be a function $\textstyle{x \mapsto 
f_K(x_1,\cdots,x_d)=\max\{ (1/d)\sum_{i=1}^d x_i -K, 0 \}}$, $K=60,70,\cdots,140$.

 We compute $(Q_{T/n}^{(m)})^n f(x_0)$, $x_0=(x_{0,1},\cdots,x_{0,d})=(100,\cdots,100)$ with $m=2, 3$ and $n=2^0,2^1,2^2$ by the machine learning-based method for the second and the third-order weak approximation (WA2.0 and WA3.0, for short) where the batch size and the iteration steps are taken as $M=256$, $1024$ and $J=4000$, respectively. The learning rate is taken as $\gamma\equiv \gamma_j$, $j\leq J$, with $\gamma=5 \times 10^{-1}{\bf 1}_{[0,600]} +10^{-2}{\bf 1}_{(600,1200]} +10^{-3}{\bf 1}_{(1200,4000]}$ for $K<90$, $\gamma=10^{-1}{\bf 1}_{[0,600]} +10^{-2}{\bf 1}_{(600,1200]} +10^{-3}{\bf 1}_{(1200,4000]} $ for $90\leq K\leq 100$ and $\gamma=10^{-2}{\bf 1}_{[0,600]} +10^{-3}{\bf 1}_{(600,1200]} +10^{-4}{\bf 1}_{(1200,4000]}$ for $100<K$. The reference value of $P_T f_K(x_0)$ is computed by Monte Carlo simulation with the number of paths $10^8$ ($\mathrm{Ref}(K)$) by using explicit solution of $X^x$ obtained by It\^o formula.  Also, for comparison, we compute numerical values obtained by the Euler-Maruyama scheme by the machine learning-based method (EM ML) with the same batch size and the same iteration steps in WA2.0 and WA3.0 where the number of discretization is taken as $n=2^0,2^1,\cdots,2^{11}$, which corresponds to machine learning-based computation for $(Q_{T/n}^{(m)})^n f(x_0)$ with $m=1$. Furthermore, as a reference, numerical values of the standard method, the Euler-Maruyama scheme with the standard Monte Carlo simulation (EM MC) are shown, where the number of discretization and the number of paths are taken as $n=2^0,2^1,\cdots,2^{11}$ and $10^5, 10^6$.  We compute the mean of each WA2.0, WA3.0 and EM ML (EM MC) (and runtime) of 25 independent trials and plot the absolute error of each method. More precisely, we test that how the maximum absolute errors ${\cal E}^{\mathrm{WA2.0}}(n)$, ${\cal E}^{\mathrm{WA3.0}}(n)$ and ${\cal E}^{\mathrm{EM \ ML}}(n)$ (${\cal E}^{\mathrm{EM \ MC}}(n)$) decrease with respect to time-steps $n$, where ${\cal E}^{\mathrm{WA2.0}}(n)$, ${\cal E}^{\mathrm{WA3.0}}(n)$ and ${\cal E}^{\mathrm{EM \ ML}}(n)$ (${\cal E}^{\mathrm{EM \ MC}}(n)$) are given by ${\cal E}^{\mathrm{WA2.0}}(n)=|\mathrm{Mean}^{\mathrm{WA2.0}(k,n)}-\mathrm{Ref}(k)|$ with $k=\mathrm{argmax}_K |\mathrm{Mean}^{\mathrm{WA2.0}(K,2^0)}-\mathrm{Ref}(K)|$, ${\cal E}^{\mathrm{WA3.0}}(n)=|\mathrm{Mean}^{\mathrm{WA3.0} (k,n)}-\mathrm{Ref}(k)|$ with $k=\mathrm{argmax}_K |\mathrm{Mean}^{\mathrm{WA3.0}(K, 2^0)}-\mathrm{Ref}(K)|$ and ${\cal E}^{\mathrm{EM \ ML}}(n)=|\mathrm{Mean}^\mathrm{EM \ ML}(k,n)-\mathrm{Ref}(k)|$ (${\cal E}^{\mathrm{EM \ MC}}(n)=|\mathrm{Mean}^\mathrm{EM \ MC}(k,n)-\mathrm{Ref}(k)|$) with $k=\mathrm{argmax}_K |\mathrm{Mean}^{\mathrm{EM \ ML}(K,2^0)}-\mathrm{Ref}(K)|$ ($k=\mathrm{argmax}_K |\mathrm{Mean}^{\mathrm{EM \ MC}(K,2^0)}-\mathrm{Ref}(K)|$). Here, the labels $\mathrm{Mean}^{\mathrm{WA2.0}(K,n)}$, $\mathrm{Mean}^{\mathrm{WA3.0}(K,n)}$ and $\mathrm{Mean}^{\mathrm{EM \ ML}(K,n)}$ ($\mathrm{Mean}^{\mathrm{EM \ MC}(K,n)}$) are the mean of 25 independent trials of $\mathrm{WA2.0}(K,n)$, $\mathrm{WA3.0}(K,n)$ and $\mathrm{EM \ ML}(K,n)$ ($\mathrm{EM \ MC}(K,n)$), where $\mathrm{WA2.0}(K,n)$, $\mathrm{WA3.0}(K,n)$ and $\mathrm{EM \ ML}(K,n)$ ($\mathrm{EM \ MC}(K,n)$) represent the approximate values for  $E[f_K(X_T^{x_0})]$ computed by WA2.0, WA3.0 and EM ML (EM MC) with the number of time steps $n$. 
 
The results are plotted in Figure 1. In the figure, the labels ``WAx.0 (SGD batch size xxx, train steps yyy)" and ``EM (SGD batch size xxx, train steps yyy)" represent the 
 x-th order weak approximation and EM ML respectively, computed with the batch size xxx used in each experiment for the machine learning-based method. Also, the labels such as ``EM (MC xxx paths)" are used for EM MC with the number of paths xxx. By the figure, we see that WA2.0 and WA3.0 with the batch size $M=1024$ and EM with MC paths $10^6$ give stable results although WA2.0 and WA3.0 with batch size $M=256$ and EM with MC paths $10^5$ provide poor estimates. Also, it can be checked that WA2.0 and WA3.0 with batch size $M=1024$ attain second and third-order convergence and give enough accuracy with smaller number of 
 time steps with short runtimes. Actually, WA2.0 and WA3.0 with $n=2^2$ and batch size $M=1024$ are computed with $7.29$s and $8.82$s, respectively while EM ML with $n=2^{11}$ and batch size $M=1024$ takes $1635.98$s (EM MC with $n=2^{11}$ and paths 
 $10^6$ (using matrix computation with NumPy) takes $990.15$s) to be the same error level.  

\begin{figure}[H]
\centering  \includegraphics[scale=0.45]{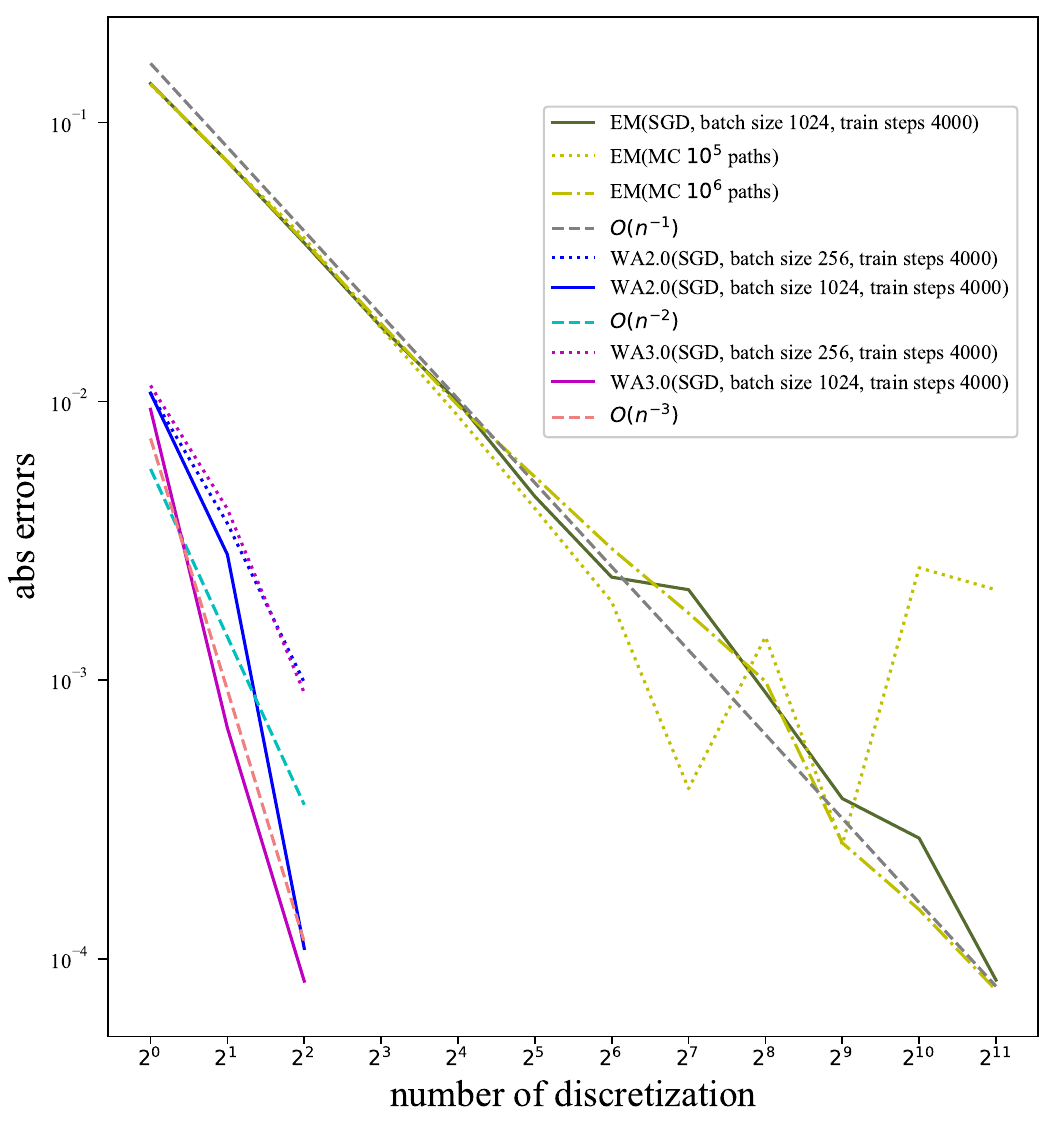}
\caption{Rate of convergence of weak approximation in 10-dimensional SDE ($d=10$)
}
\end{figure}

\subsection{100-dimension case}
Next, we show numerical results for the $100$-dimension case. Let $d=100$, $\sigma=0.2$, $T=1.0$, $\textstyle{f=f_K:x \mapsto f_K(x_1,\cdots,x_d)=\max\{ \max\{ x_1-K, 0 \}, \cdots,\max\{ x_d -K, 0 \} \} }$, $K=60,70,\cdots,140$. We compute $(Q_{T/n}^{(m)})^n f(x_0)$, $x=(x_{0,1},\cdots,x_{0,d})=(100,\cdots,100)$ with $m=2, 3$ and $n=2^i$ up to $i=2$ or $3$ using WA2.0 and WA3.0 where the batch size and the iteration steps are taken as $M=1024$ and $J=2000,4000$, respectively. The learning rate is taken as $\gamma\equiv \gamma_j$, $j\leq J$, with $\gamma=5 \times 10^{-1}{\bf 1}_{[0,600]} +5\times 10^{-2}{\bf 1}_{(600,1200]} +5\times 10^{-3}{\bf 1}_{(1200,4000]}$. The benchmark value of $P_T f(x_0)$ is computed by MC with the number of paths $10^8$ using explicit solution of $X^{x_0}$ obtained by It\^o formula.  Also, for comparison, we compute the numerical values obtained by the Euler-Maruyama scheme (EM ML and EM MC) where the number of discretization is taken as $n=2^0,2^1,\cdots,2^{9}$.  We perform same experiments for this case.

 The results are plotted in Figure 2.  We again check that WA2.0 and WA3.0 attain second and third-order convergence. The schemes WA2.0 and WA3.0 give enough accuracy with smaller number of time steps with short runtimes. Actually, WA3.0 with $n=2^2$, batch size $M=1024$ with iteration steps $J=2000$ and $J=4000$ are computed with runtime $30.33$s ($J=2000$) and $60.18$s ($J=4000$) respectively, while EM ML with $n=2^{9}$, batch size $M=1024$   with iteration steps $J=2000$ and $J=4000$ take $2328.77$s ($J=2000$) and $4702.91$s ($J=4000$) (also EM MC with $n=2^{9}$ and paths $10^6$ (using matrix computation with NumPy) takes $2483.44$s) to be the same error level.

\begin{figure}[H]
\centering  \includegraphics[scale=0.45]{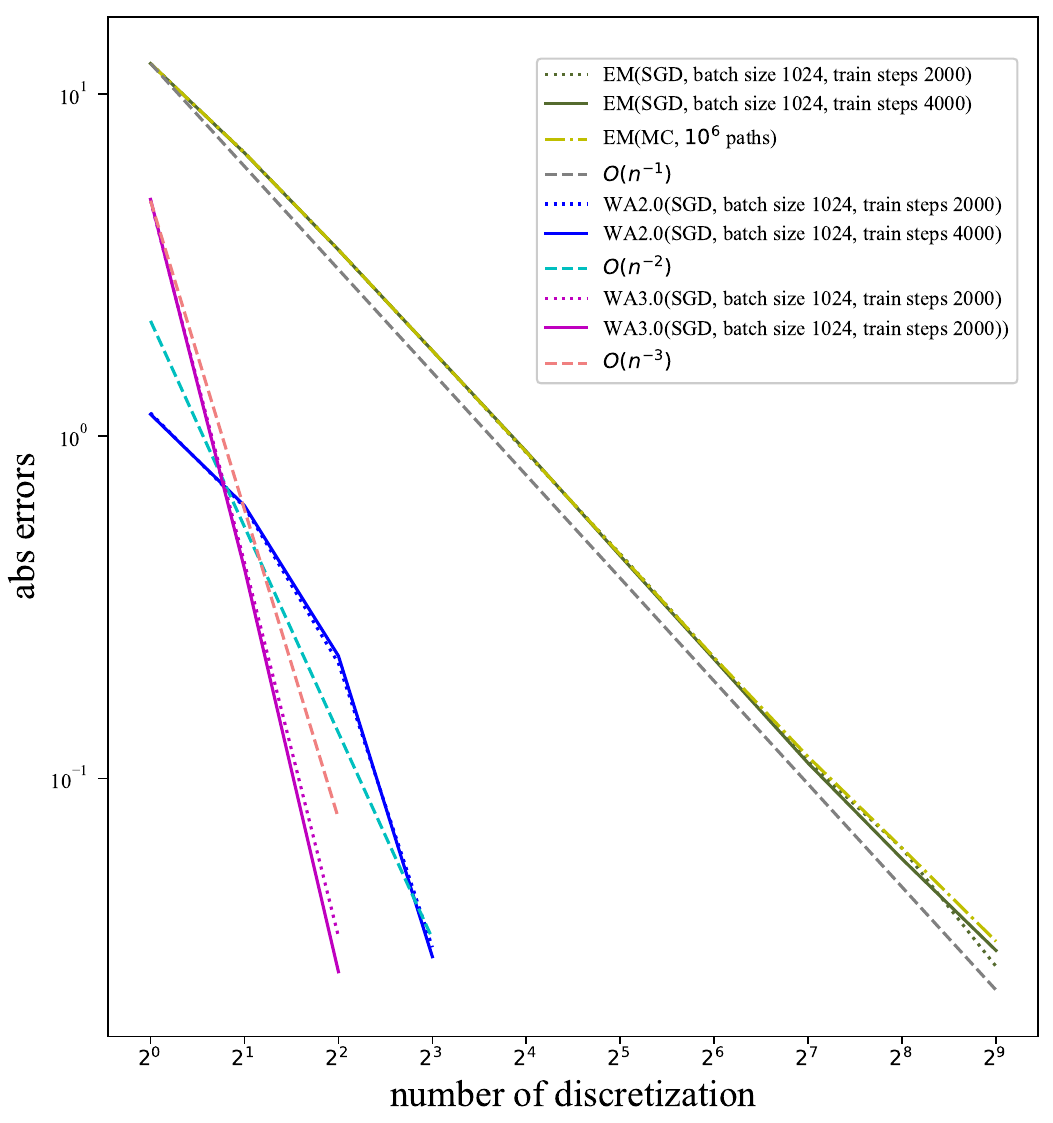}
\caption{Rate of convergence of weak approximation in 100-dimensional SDE ($d=100$)
}
\end{figure}

Through the experiments, we have checked that the machine learning-based high-order weak approximation can be a new computation method for high-dimensional integrals.

\section{Conclusion and future works}
In the paper, we introduced a new machine learning-based high-order weak approximation method for solving high-dimensional Kolmogorov PDEs. The target high-dimensional integral is efficiently computed through stochastic weighted minimization and stochastic gradient descent method without suffering from the curse of dimensionality. The proposed method may be an alternative for existing methods such as finite difference method, finite element method and Monte Carlo method. 

We need more experiments for optimal choice of parameters in the stochastic gradient descent approximation. Also, improved stochastic gradient descent algorithm such as \cite{LJHCLGH} should be applied to the proposed method of the paper. 

We believe that  ``deep learning" enables us to give wide applications in which the ``function $(Q_{T/n}^{(m)})^n f$" is required. In such cases, we need to combine the proposed method with deep learning, while for a fixed $x_0$, it is enough to use the proposed machine learning-based method for computing the value of $(Q_{T/n}^{(m)})^n f(x_0)$. As future work, we will study whether the method is applicable to solve high-dimensional nonlinear PDEs (the problem as in \cite{BBCJN}) with short runtimes by means of deep learning.   

\section*{Acknowledgements}
The second author thanks Prof. Masaaki Fukasawa, Prof. Jun Sekine and Prof. Akihiko Takahashi for providing useful comments on the method of this paper. This work is supported by JSPS KAKENHI (Grant Number 19K13736) and JST PRESTO (Grant Number JPMJPR2029), Japan.

\end{document}

%% file: tymacro_3.tex

\def\bn{{\bf n}}
\def\A{{\bf A}}
\def\B{{\bf B}}
\def\C{{\bf C}}
\def\D{{\bf D}}
\def\E{{\bf E}}
\def\F{{\bf F}}
\def\G{{\bf G}}
\def\H{{\bf H}}
\def\I{{\bf I}}
\def\J{{\bf J}}
\def\K{{\bf K}}
\def\L{{\bf L}}
\def\M{{\bf M}}
\def\N{{\bf N}}
\def\O{{\bf O}}
\def\P{{\bf P}}
\def\Q{{\bf Q}}
\def\R{{\bf R}}
\def\S{{\bf S}}
\def\T{{\bf T}}
\def\U{{\bf U}}
\def\V{{\bf V}}
\def\W{{\bf W}}
\def\X{{\bf X}}
\def\Y{{\bf Y}}
\def\Z{{\bf Z}}
\def\cala{{\cal A}}
\def\calb{{\cal B}}
\def\calc{{\cal C}}
\def\cald{{\cal D}}
\def\cale{{\cal E}}
\def\calf{{\cal F}}
\def\calg{{\cal G}}
\def\calh{{\cal H}}
\def\cali{{\cal I}}
\def\calj{{\cal J}}
\def\calk{{\cal K}}
\def\call{{\cal L}}
\def\calm{{\cal M}}
\def\caln{{\cal N}}
\def\calo{{\cal O}}
\def\calp{{\cal P}}
\def\calq{{\cal Q}}
\def\calr{{\cal R}}
\def\cals{{\cal S}}
\def\calt{{\cal T}}
\def\calu{{\cal U}}
\def\calv{{\cal V}}
\def\calw{{\cal W}}
\def\calx{{\cal X}}
\def\caly{{\cal Y}}
\def\calz{{\cal Z}}
%
\def\sskip{\hspace{0.5cm}}
\def\simleq{ \raisebox{-.7ex}{\em $\stackrel{{\textstyle <}}{\sim}$} }
\def\leqsim{ \raisebox{-.7ex}{\em $\stackrel{{\textstyle <}}{\sim}$} }
\def\ep{\epsilon}
\def\half{\frac{1}{2}}
\def\iku{\rightarrow}
\def\Iku{\Rightarrow}
\def\ikup{\rightarrow^{p}}
\def\inclusion{\hookrightarrow}
\def\cadlag{c\`adl\`ag\ }
\def\up{\uparrow}
\def\down{\downarrow}
\def\doti{\Leftrightarrow}
\def\douti{\Leftrightarrow}
\def\dochi{\Leftrightarrow}
\def\douchi{\Leftrightarrow}%
\def\yy{\\ && \nonumber \\}
\def\y{\vspace*{3mm}\\}
\def\nn{\nonumber}
\def\be{\begin{equation}}
\def\ee{\end{equation}}
\def\bea{\begin{eqnarray}}
\def\eea{\end{eqnarray}}
\def\beas{\begin{eqnarray*}}
\def\eeas{\end{eqnarray*}}
%
\def\hd{\hat{D}}
\def\hv{\hat{V}}
\def\hsd{{\hat{d}}}
\def\hx{\hat{X}}
\def\hsx{\hat{x}}
\def\bsx{\bar{x}}
\def\bsd{{\bar{d}}}
\def\bx{\bar{X}}
\def\ba{\bar{A}}
\def\bb{\bar{B}}
\def\bc{\bar{C}}
\def\bv{\bar{V}}
\def\balpha{\bar{\alpha}}
\def\bbalpha{\bar{\bar{\alpha}}}
\def\combi{\l(\begin{array}{c}\alpha\\ \beta \end{array}\r)}
\def\f{^{(1)}}
\def\s{^{(2)}}
\def\ss{^{(2)*}}
\def\l{\left}
\def\r{\right}
\def\a{\alpha}
\def\b{\beta}
\def\L{\Lambda}